# A three point formula for finding roots of equations by the method of least squares

Ababu Teklemariam Tiruneh[1] ; William N. Ndlela[1] ; Stanley J. Nkambule[1]

[1] Lecturer, Department of Environmental Health Science. University of Swaziland.
P.O.Box 369, Mbabane H100, Swaziland.  Email: ababute@yahoo.com

**Abstract**

A new method of root finding is formulated that uses a numerical iterative process involving three points. A given function Y= f(x) whose root(s) are desired is fitted and approximated by a polynomial function curve of the form y= a(x-b)$^N$ and passing through three equi-spaced points using the method of least squares.  Successive iterations using the same procedure of curve fitting is used to locate the root within a given level of tolerance.  The power N of the curve suitable for a given function form can be appropriately varied at each step of the iteration to give a faster rate of convergence and avoid cases where oscillation, divergence or off shooting to an invalid domain may be encountered. An estimate of the rate of convergence is provided. It is shown that the method has a quadratic convergence similar to that of Newton's method. Examples are provided showing the procedure as well as comparison of the rate of convergence with the secant and Newton's methods. The method does not require evaluation of function derivatives.

## 1. Introduction

Finding the roots of equations through numerical iterative procedure is an important step in the solution of many science and engineering problems.  Beginning with the classical Newton method, several methods for finding roots of equations have been proposed each of which has its own advantages and limitations. Newton's method of root finding is based on the iterative formula:

$$x_{k+1} = x_k - \frac{y(x_k)}{y'(x_k)}$$

Newton's method has a  quadratic convergence and requires a derivative of the function for each step of the iteration. When the derivative evaluated is zero, Newton method fails. For low values of the derivative the Newton iteration offshoots away from the current point of iteration. The convergence of Newton method can be slow near roots of multiplicity although modifications can be made to increase the rate of convergence [1].

Acceleration of  Newton's method with higher order convergence have been proposed that require also evaluation of a function and its derivatives. For example  a third order convergence method by S. Weeraksoon and T.G. Fernando [2] requires evaluation of one function and two first derivatives. A fourth order iterative method, accordintg to J.F. Traub [3]  also requires evaluation of one function and two derivative evaluations. Sanchez and Barrero [4] gave a compositing of function evaluation at a point and its derivative to improve the convergence of Newton's method from 2 to 4. Recently other methods of fifth, sixth, seventh and higher order convergence have been proposed [5-11]. In all of such methods evaluation of function and its derivatives are necessary.

The secant method does not require evaluation of derivatives. However, the rate of convergence is about 1.618. Muller's method is an extension of the secant method to a quadratic polynomial [12]. It requires three functional evaluations to start with but continues with one function evaluation afterwards. The method does not require derivatives and the rate of convergence is about 1.84. However, Muller's method can converge to a complex root from an initial real number [13].

## 2. Method development

For a given function of the form $Y = f(x)$, three starting points separated by an equi-spaced horizontal distance of $\delta$ are chosen. The points pass through the given function $Y = f(x)$. A single root polynomial function of the general form $Y = a(x-b)^N$ is fitted to the given points using the method of least squares. $N$ is the power of the polynomial which is generally a real number and b is the root of the polynomial which serves to approximate the root of the given function $y = f(x)$ at any given step of the iteration process. Figure 1 shows the three different possible curves that can be fitted to a given function using the three points.

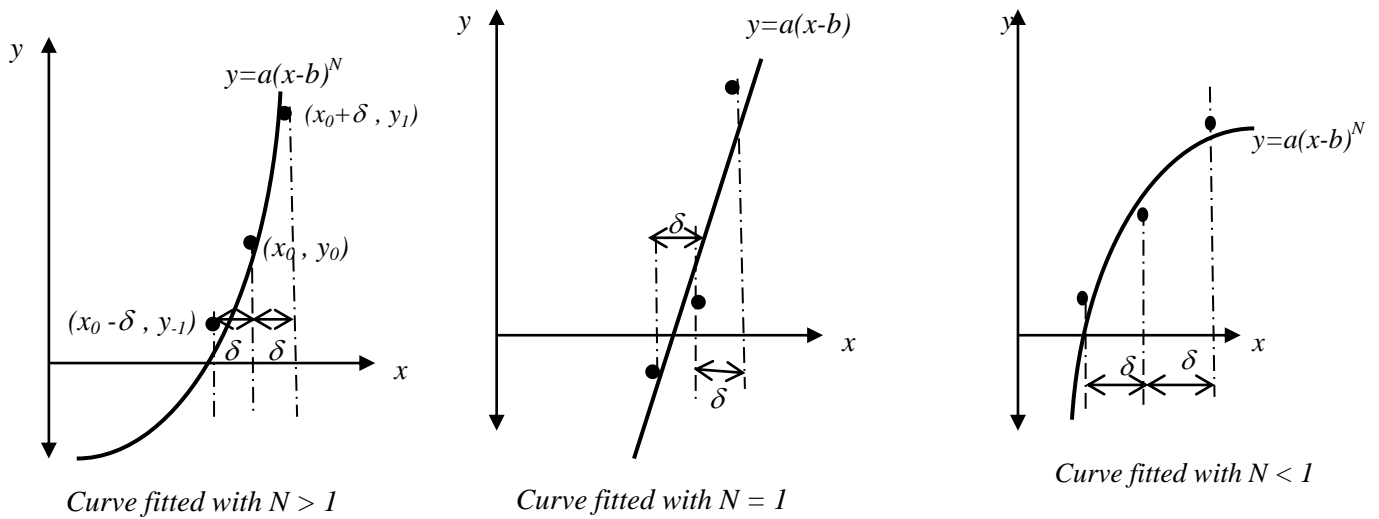

*Curve fitted with N > 1*    *Curve fitted with N = 1*    *Curve fitted with N < 1*

Figure 1. Different types of curves that can be fitted to the three points using the method of least squares.

Depending on the behavior of the function $Y = f(x)$ to be approximated, the power of the polynomial N, where N is generally a real number, can be chosen. Figure 1 above shows N can take values greater than 1, equal to one or can be less than one  The constants *a* and *b* are determined by applying the method of least squares by minimizing the sum of the squares of the errors in y values over the three points namely $(x_0 -\delta , y_{-1})$ $(x_0 , y_0)$ and $(x_0 +\delta , y_1)$.

$$\sum_{i=-1}^{i=1} e_i^2 = \sum_{i=-1}^{i=1} [y_i - a(x_i - b)^N]^2 \tag{1}$$

Differentiating Equation 1 with respect to *a* and setting the resulting expression to zero will give the following expression for the constant *a*:

$$a = \frac{\sum_{i=-1}^{i=1} y_i (x_i - b)^N}{\sum_{i=-1}^{i=1} (x_i - b)^{2N}} \tag{2}$$

Differentiation of the sum of squares of the errors (Equation 1) again with respect to the constant *b* and setting the expression to zero will also give the following equivalent expression for *a*

$$a = \frac{\sum_{i=-1}^{i=1} y_i (x_i - b)^{N-1}}{\sum_{i=-1}^{i=1} (x_i - b)^{2N-1}} \tag{3}$$

The constant *a* is not desired for root approximation. The root approximation of the polynomial curve $y = a(x-b)^N$ is the constant *b*. Therefore, the above two equations (equations 2 and 3) are equated to eliminate *a*, resulting in the following expression:

$$\left(\sum_{i=-1}^{i=1} y_i(x_i - b)^N\right)\left(\sum_{i=-1}^{i=1} (x_i - b)^{2N-1}\right) = \left(\sum_{i=-1}^{i=1} y_i(x_i - b)^{N-1}\right)\left(\sum_{i=-1}^{i=1} (x_i - b)^{2N}\right) \tag{4}$$

Choosing the three points that are equi-spaced separated by a horizontal distance of $\delta$ will result in simplified expression for the root b. Therefore, the three points $(x_{-1}, y_{-1})$, $(x_0, y_0)$ and $(x_1, y_1)$ are replaced by $(x_0-\delta , y_{-1})$ $(x_0 , y_0)$ and $(x_0+\delta , y_1)$ respectively.

The bracketed expressions in the equation above are each evaluated by making use of binomial expansion of the terms involving $x_0$, $\delta$ and b raised to the various powers of N. In addition, for small values of $\delta$, the terms containing $\delta^3$ and higher orders are discarded. The resulting expressions are the following:

$$\sum_{i=-1}^{i=1} y_i(x_i - b)^N = \left(\sum_{i=-1}^{i=1} y_i\right)(x_0 - b)^N + N\delta(y_1 - y_{-1})(x_0 - b)^{N-1} + \left(\frac{N(N - 1)\delta^2}{2}\right)(y_1 + y_{-1})(x_0 - b)^{N-2} \tag{5}$$

$$\sum_{i=-1}^{i=1} (x_i - b)^{2N-1} = 3(x_0 - b)^{2N-1} + 2(N - 1)(2N - 1)\delta^2 (x_0 - b)^{2N-3} \tag{6}$$

$$\sum_{i=-1}^{i=1} y_i(x_i - b)^{N-1} = \left(\sum_{i=-1}^{i=1} y_i\right)(x_0 - b)^{N-1} + (N-1)\delta(y_1 - y_{-1})(x_0 - b)^{N-2} + \left(\frac{(N-1)(N-2)\delta^2}{2}\right)(y_1 + y_{-1})(x_0 - b)^{N-3} \quad (7)$$

$$\sum_{i=-1}^{i=1}(x_i - b)^{2N} = 3(x_0 - b)^{2N} + 2(N)(2N - 1)\delta^2 (x_0 - b)^{2N-2} \quad (8)$$

Substituting the above expressions for the bracketed products of Equation 4 and again discarding the terms containing $\delta^3$ and higher order (for small values of $\delta$), gives the following expression:

$$[3N(y_1 - y_{-1})(x_0 - b)^{3N-2}]\delta - [3(N-1)(y_1 - y_{-1})(x_0 - b)^{3N-2}]\delta =$$

$$\left[2N(2N-1)\left(\sum_{i=-1}^{i=1} y_i\right)((x_0 - b)^{3N-3})\right]\delta^2 - \left[2(N-1)(2N-1)\left(\sum_{i=-1}^{i=1} y_i\right)((x_0 - b)^{3N-3})\right]\delta^2 +$$

$$\left[\left(\frac{3(N-1)(N-2)}{2}\right)(y_1 + y_{-1})(x_0 - b)^{3N-3}\right]\delta^2 - \left[\left(\frac{3(N)(N-2)}{2}\right)(y_1 + y_{-1})(x_0 - b)^{3N-3}\right]\delta^2 \quad (9)$$

Since $\delta$ is common to all the expressions, it is factored out from all the terms, and, further simplification leads to:

$$3(y_1 - y_{-1})(x_0 - b) = \left[(4N - 2)\left(\sum_{i=-1}^{i=1} y_i\right)\right]\delta + [(-3N + 3)(y_1 + y_{-1})]\delta \quad (10)$$

Solving for the variable $b$, which is the approximation to the root at a given iteration, gives:

$$b = x_0 - N\left[\frac{\left(\frac{(N+1)y_{-1} + (4N-2)y_0 + (N+1)y_1}{6N}\right)}{\left(\frac{y_1 - y_{-1}}{2\delta}\right)}\right] \quad (11)$$

In terms of the iteration process the estimate of the root at the $(k+1)^{th}$ iteration is evaluated from functional values of the $k^{th}$ iteration, the above expression can be written as:

$$x_{k+1} = x_k - N\left[\frac{\left(\frac{(N+1)y_{k-\delta} + (4N-2)y_k + (N+1)y_{k+\delta}}{6N}\right)}{\left(\frac{y_{k+\delta} - y_{k-\delta}}{2\delta}\right)}\right] \quad (12)$$

It is interesting to mention the similarity with Newton's expression for approximation of roots. The numerator in the bracket when multiplied is the weighted average of y values with the central point having a weight of 4N-2 while the end points each are weighed by N+1. The denominator in the square bracket is the central difference approximation to the derivative for the central point ($x_k$, $y_k$). The N value outside the brackets represents the 'acceleration' factor as in the Newton's method whereby the iteration accelerates when the factor N is applied to the Newton method of root finding, for example, for roots of polynomials with root multiplicity of N.

The equivalent expressions for N = 1, 2 and 3 are given as follows:

$$\text{For } N = 1 \quad ; \quad x_{k+1} = x_k - 1 \left[ \frac{\left(\frac{y_{k-\delta} + y_k + y_{k+\delta}}{3}\right)}{\left(\frac{y_{k+\delta} - y_{k-\delta}}{2\delta}\right)} \right] \quad (13)$$

$$\text{For } N = 2 \quad ; \quad x_{k+1} = x_k - 2 \left[ \frac{\left(\frac{y_{k-\delta} + 2y_k + y_{k+\delta}}{4}\right)}{\left(\frac{y_{k+\delta} - y_{k-\delta}}{2\delta}\right)} \right] \quad (14)$$

$$\text{For } N = 3 \quad ; \quad x_{k+1} = x_k - 3 \left[ \frac{\left(\frac{2y_{k-\delta} + 5y_k + 2y_{k+\delta}}{9}\right)}{\left(\frac{y_{k+\delta} - y_{k-\delta}}{2\delta}\right)} \right] \quad (15)$$

**Estimation of the power of the polynomial N at each step of the iteration**

It is possible to vary the power of the polynomial N in the equation $y = a(x-b)^N$ in each iteration step which means different curves can be fitted depending on the curve defined by the three points. The estimated value of N to be used in the iteration formula will be derived from the y-values ($y_{k-\delta}$, $y_k$, $y_{k+\delta}$) without involving any of the derivatives. However, for the purpose of derivation of N, the derivatives will be used which will be eventually replaced by the finite difference form approximations.

From the equation $y = a(x-b)^N$, the first derivative dy/dx and $2^{nd}$ derivative $d^2y/dx^2$ are given by;

$$\frac{dy}{dx} = Na(x-b)^{N-1} \quad (16)$$

$$\frac{d^2y}{dx^2} = N(N-1)a(x-b)^{N-2} \quad (17)$$

From the expression of y, dy/dx and $d^2y/dx^2$ above the following two equations are obtained:

$$\frac{y}{dy/dx} = \frac{a(x-b)^N}{Na(x-b)^{N-1}} = \frac{(x-b)}{N} \quad (18)$$

$$\frac{dy/dx}{\frac{d^2y}{dx^2}} = \frac{Na(x-b)^{N-1}}{N(N-1)a(x-b)^{N-2}} = \frac{(x-b)}{N-1} \tag{19}$$

Eliminating (x-b) from Equation 18 and equation 19 and solving for the power N gives

$$N = \frac{(dy/dx)^2}{(dy/dx)^2 - y\left(\frac{d^2y}{dx^2}\right)} \tag{20}$$

Replacing the derivatives by the finite difference approximations involving the three equidistant points ($y_{k-\delta}$, $y_k$, $y_{k+\delta}$):

$$N_k = \frac{\left(\frac{y_{k+\delta} - y_{k-\delta}}{2\delta}\right)^2}{\left(\frac{y_{k+\delta} - y_{k-\delta}}{2\delta}\right)^2 - y_k\left(\frac{y_{k-\delta} - 2y_k + y_{k+\delta}}{\delta^2}\right)} \tag{21}$$

**Proof of quadratic convergence**

Recalling the root approximation formula in an iteration form involving the $k^{th}$ and $(k+1)^{th}$ iterations (i.e. Eq. 12):

$$x_{k+1} = x_k - N\left[\frac{\left(\frac{(N+1)y_{k-\delta} + (4N-2)y_k + (N+1)y_{k+\delta}}{6N}\right)}{\left(\frac{y_{k+\delta} - y_{k-\delta}}{2\delta}\right)}\right]$$

Expanding $y_{k-\delta}$ and $y_{k+\delta}$ about $y_k$ using Taylor series expansion;

$$y_{k+\delta} = y_k + \delta y'_k + \frac{\delta^2 y''_k}{2} + \frac{\delta^3 y'''_k}{6} + O(\delta)^4 \tag{22}$$

$$y_{k-\delta} = y_k - \delta y'_k + \frac{\delta^2 y''_k}{2} - \frac{\delta^3 y'''_k}{6} + O(\delta)^4 \tag{23}$$

Inserting the above expression in the numerator of the iteration formula yields;

$$\left(\frac{(N+1)y_{k-\delta} + (4N-2)y_k + (N+1)y_{k+\delta}}{6N}\right) = y_k + \left(\frac{N+1}{6N}\right)\delta^2 y''_k + O(\delta)^4 \tag{24}$$

Similarly the denominator will, after substitution of the Taylor series expression, reduces to;

$$\left(\frac{y_{k+\delta} - y_{k-\delta}}{2\delta}\right) = y'_k + \frac{\delta^2 y'''_k}{6} + O(\delta)^3 \tag{25}$$

Assuming the expression $\frac{\delta^2 y_k'''}{6} + O(\delta)^3$ to be small compared to y'$_k$ and neglecting this term will give:

$$x_{k+1} = x_k - N\left(\frac{y_k + \left(\frac{N+1}{6N}\right)\delta^2 y''_k}{y'_k}\right) + O(\delta)^4 \qquad (26)$$

Defining the error $E_k$ at the k$^{th}$ iteration as $E_k = x_k - r$ where r is the root. Also $E_{k+1} = x_{k+1} - r$
Substituting $E_k + r$ for $x_k$ and $E_{k+1} + r$ for $x_{k+1}$ yields;

$$E_{k+1} = E_k - N\left(\frac{y_k + \left(\frac{N+1}{6N}\right)\delta^2 y''_k}{y'_k}\right) + O(\delta)^4 \qquad (27)$$

The above expression will be worked out further for two different cases. The first case is for N=1 and the second for N is any number different from 1 and providing root of multiplicity N.

For the first case, N=1 ;

$$E_{k+1} = E_k - \left(\frac{y_k + \left(\frac{1}{3}\right)\delta^2 y''_k}{y'_k}\right) + O(\delta)^4 \qquad (28)$$

Expanding y$_k$ about the root x= r , using Taylor series expansion where $x_k = E_k + r$ ;

$$y_k = (y_r = 0) + E_k y'_r + \frac{E_k^2 y''_r}{2} + O(E_k)^3 \qquad (29)$$

Similarly expanding y'$_k$ about r $_{gives}$;

$$y'_k = y'_r + E_k y''_r + \frac{E_k^2 y''_r}{2} + O(E_k)^3 \qquad (30)$$

Again assuming the terms $E_k y''_r + \frac{E_k^2 y''_r}{2} + O(E_k)^3$ to be small compared to y'$_r$

$$y'_k = y'_r \qquad (31)$$

Similarly expansion of y"$_k$ about r gives;

$$y''_k = y''_r + E_k y'''_r + \frac{E_k^2 y_r^{iv}}{2} + O(E_k)^3 \qquad (32)$$

Substituting for y$_k$ , y'$_k$ and y"$_k$ gives;

$$E_{k+1} = E_k - \left[\frac{E_k y'_r + \frac{E_k^2 y''_r}{2} + O(E_k)^3 + \left(\frac{1}{3}\right)\delta^2 \left(y''_r + E_k y'''_r + \frac{E_k^2 y_r^{iv}}{2} + O(E_k)^3\right)}{y'_r}\right] + O(\delta)^4 \quad (33)$$

Reducing further gives;

$$E_{k+1} = -\frac{E_k^2 y''_r}{2y'_r} - \delta^2 \frac{y''_r}{3y'_r} - \delta^2 E_k \frac{y'''_r}{3y'_r} - \delta^2 E_k^2 \frac{y_r^{iv}}{6y'_r} + O(E_k)^3 + O(\delta)^4 \quad (34)$$

The above expression results in convergence which is a function of $\delta^2, E_k^2, \delta^2 E_k,$ and $\delta^2 E_k^2$. The $\delta$ value is set as the square of the difference in x values of the previous successive iterations multiplied by a factor $\beta$ which is given a value less than or equal to one.

$$\delta = \beta (E_k - E_{k-1})^2 \quad (35)$$

The value of $\delta^2$ will therefore be ;

$$\delta^2 = \beta^2 (E_k - E_{k-1})^4$$

So the error series will take the form:

$$E_{k+1} = -\frac{E_k^2 y''_r}{2y'_r} - (E_k - E_{k-1})^4 \left(\frac{\beta^2 y''_r}{3y'_r}\right) - (E_k - E_{k-1})^4 E_k \left(\frac{\beta^2 y'''_r}{3y'_r}\right)$$
$$- (E_k - E_{k-1})^4 E_k^2 \left(\frac{\beta^2 y_r^{iv}}{6y'_r}\right) + O(E_k)^3 + O(\delta)^4 \quad (36)$$

It will now be shown that the $\pm (E_k - E_{k-1})^4$ term is quadratically convergent.

Assuming the $(E_k - E_{k-1})^4$ is the dominant term in the above expression which means

$$E_{k+1} = \pm(E_k - E_{k-1})^4 \geq E_k^2 \quad (37)$$

$$E_{k+1} = E_k^n = \pm(E_k - E_{k-1})^4 \quad (38)$$

For the case positive case, i.e.,
$$E_{k+1} = E_k^n = (E_k - E_{k-1})^4$$

$$E_k^{n/4} = E_k^{1/n} - E_k \quad (39)$$

The right hand term of the above expression is evaluated for small values of $E_k$ and for the following conditions;

For n = 1 ; $E_k^{1/4} = E_k - E_k = 0$ is not a valid expression.
For n < 1 ; $E_k \gg E_k^{1/n}$ so that $E_k^{n/4} = -E_k$ is also not valid expression.
For n > 1   $E_k^{1/n} \gg E_k$ so that;

$$E_k^{n/4} = E_k^n$$

$$\frac{n}{4} = \frac{1}{n}$$

$$n^2 = 4 \quad \text{or} \quad n = 2$$

Therefore, for positive $(E_k - E_{k-1})^4$ term, $E_{k+1} \, \alpha \, E_k^2$

For the negative $(E_k - E_{k-1})^4$ term :

$$E_{k+1} = E_k^n = -(E_k - E_{k-1})^4 \tag{40}$$

$$E_k^{n/4} + E_k^{1/n} - E_k = 0 \tag{41}$$

Let a function f(n) be defined so that:

$$f(n) = E_k^{n/4} + E_k^{1/n} - E_k \tag{42}$$

It is possible to show that for all n ≥ 0 the function f(n) is always positive or always negative depending on the sign of $E_k$. To show this the following ranges are considered:

For $0 \leq n \leq 1$  $E_k^{n/4}$ is the dominant term so that $f(n) = E_k^{n/4}$
For $1 < n < 4$ ,  $E_k^{n/4} + E_k^{1/n}$ are both dominant so that $f(n) = E_k^{n/4} + E_k^{1/n}$
For $n \geq 4$ ,  $E_k^{1/n}$ is the dominant term so that $f(n) = E_k^{1/n}$

Therefore, the equation f(n) = 0 represents the minimum in the case of positive f(n) and the maximum in the case of negative f(n) values. The value of n is then determined for maximum or minimum case by setting its derivative to zero, i.e.,

$$\frac{df}{dn} = 0 \tag{43}$$

$$\frac{df}{dn} = \frac{1}{4} E_k^{n/4} \ln(E_k) - \frac{1}{n^2} E_k^{\frac{1}{n}} \ln(E_k) = 0 \tag{44}$$

$$E_k^{n/4} = \frac{4}{n^2} E_k^{1/n}$$

Equating the powers of $E_k$ results in:

$$\frac{n}{4} = \frac{1}{n}$$

$$n^2 = 4 \quad \text{or} \quad n = 2$$

Also the coefficient $4/n^2 = 1$ for n=2 making the expression:

$$E_k^{n/4} = \frac{4}{n^2} E_k^{1/n} \tag{45}$$

a valid expression. This proves once again that $E_{k+1} \ \alpha \ E_k^2$

A plot of the variation of f(n) for values of n between 1 and 4 for Ek= $10^{-22}$ in Figure 2 below shows the minimum value for f(n) occurs at n =2 as derived above. The function f(n) at n=2 is equal to $10^{-57}$ and is the smallest magnitude that can be attained and which occurs only by setting n=2. The function f(n) is not equal to zero as such but attains the smallest possible value (close to zero) which is made possible by setting n=2.

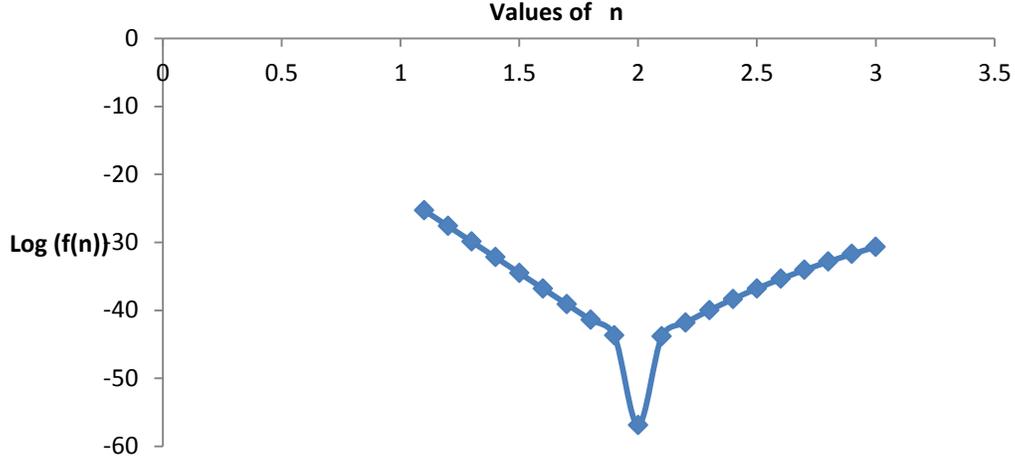

Figure 2. A plot of the error term f(n) for values of n between 1 and 4 and for $E_k = 10^{-22}$

The error series will then take the form:

$$E_{k+1} = -\frac{E_k^2 y''_r}{2y'_r} - E_k^2 \left(\frac{\beta^2 y''_r}{3y'_r}\right) - E_k^2 E_k \left(\frac{\beta^2 y'''_r}{3y'_r}\right) - E_k^2 E_k^2 \left(\frac{\beta^2 y_r^{iv}}{6y'_r}\right) + O(E_k)^3 + O(E_k^2)^4 \qquad (46)$$

$$E_{k+1} = -\frac{E_k^2 y''_r}{2y'_r} - E_k^2 \left(\frac{\beta^2 y''_r}{3y'_r}\right) - E_k^3 \left(\frac{\beta^2 y'''_r}{3y'_r}\right) - E_k^4 \left(\frac{\beta^2 y_r^{iv}}{6y'_r}\right) + O(E_k)^3 \\ + O(E_k)^8 \qquad (47)$$

$$E_{k+1} = -\left(\frac{y''_r}{2y'_r} + \frac{\beta^2 y''_r}{3y'_r}\right) E_k^2 + O(E_k)^3 \qquad (48)$$

This proves the quadratic convergence for N=1.

For N values other than one, the convergence is estimated by assuming a root of multiplicity N so that the y function is written in the form:

$$y = (x - r)^N Q(x) \qquad (49)$$

Consider the iteration formula that is in reduced form and was given by Equation 26;

$$x_{k+1} = x_k - N\left(\frac{y_k + \left(\frac{N+1}{6N}\right)\delta_k^2 y''_k}{y'_k}\right) + O(\delta_k)^4$$

The above iteration process can be written in fixed point form $x_{k+1} = g(x_k)$ by defining $g(x_k)$ such that:

$$g(x_k) = x_k - N\left(\frac{y_k + \left(\frac{N+1}{6N}\right)\delta_k^2 y''_k}{y'_k}\right) + O(\delta_k)^4 \tag{50}$$

$$g(x_k) = x_k - N\left(\frac{y_k}{y'_k}\right) - \frac{\left(\frac{N+1}{6}\right)\delta_k^2 y''_k}{y'_k} + O(\delta_k)^4 \tag{51}$$

Substituting the root of multiplicity N term $y = (x-r)^N Q(x)$ for y and the corresponding derivative of the second term of the above equation;

$$g(x_k) = x_k - N\left(\frac{N(x-r)^N Q(x)}{N(x-r)^{N-1}Q(x) + Q'(x)(x-r)^N}\right) - \frac{\left(\frac{N+1}{6}\right)\delta_k^2 y''_k}{y'_k} + O(\delta_k)^4 \tag{52}$$

It is possible to show that for the first derivative of g(x), i.e., g'(x) the first two terms cancel each other, i.e.,

$$\frac{d}{dx}\left(x_k - N\left(\frac{N(x-r)^N Q(x)}{N(x-r)^{N-1}Q(x) + Q'(x)(x-r)^N}\right)\right) = 0 \text{ at } x_k = r \tag{53}$$

Similarly at $x_k = r$, the expression;

$$\frac{d}{dx}\left(\frac{\left(\frac{N+1}{6}\right)\delta_k^2 y''_k}{y'_k} + O(\delta_k)^4\right) = 0 \tag{54}$$

holds true because $\delta_r = 0$ at the root $x_k = r$ and the derivative expressions contain the term $\delta_r = 0$ because;

$$\delta_r = \lim_{k \to \infty}(\delta_k) = \lim_{k \to \infty} \beta(E_k - E_{k-1})^2 = 0 \tag{55}$$

Therefore, g'(r) = 0

Expanding $g(x_k)$ about the root x=r using Taylor Series;

$$g(x_k) = g(r) + g'(r)E_k + \frac{g''(r)E_k^2}{2} + O(E_k)^3 \tag{56}$$

From the relation $x_{k+1} = g(x_k)$ and $r = g(r)$ and substituting $g'(r) = 0$ as shown above;

$$g(x_k) - g(r) = \frac{g''(r)E_k^2}{2} + O(E_k)^3 \tag{57}$$

$$x_{k+1} - r = E_{k+1} = g(x_k) - g(r) = \frac{g''(r)E_k^2}{2} + O(E_k)^3 \tag{58}$$

$$E_{k+1} = \frac{g''(r)E_k^2}{2} + O(E_k)^3 \tag{59}$$

Therefore, the iteration series is also quadratically convergent for N different from one.

## 3. Results and discussions

Examples of equations used to test efficiency of root finding methods are used here to evaluate the least square three-point methods and compare it case by case particularly with the Newton and secant methods. To start with, the $\delta$ value is arbitrarily set between 0 and 1 and the two points to the left and right of the central point are set as x-$\delta$ and x+$\delta$ respectively. The subsequent values of $\delta$ are set from the results of the iteration using the established formula:

$$\delta_{k+1} = \beta_{k+1}(E_k - E_{k-1})^2 = \beta_{k+1}(x_k - x_{k-1})^2 \tag{60}$$

Since the errors $E_k$ and $E_{k-1}$ are unknown the $x_k$ and $x_{k-1}$ values of the $k^{th}$ and $k-1^{th}$ iteration are used to calculate $\delta$. The $\beta_{k+1}$ value is set so that:

$$\delta_{k+1} < 1 \text{ and}$$

$$\delta_{k+1} \leq \delta_{k+1}^2$$

A simple way of reducing $\delta_{k+1}$ to satisfy the above equations is using $\beta$ value from either of the series 1, 0.1. 0.01, 0.001, etc. In most cases the use of $\beta = 1$ or $\beta = 0.1$ is adequate to satisfy the above requirements.

There are two possible options for choosing the value of the power N of the polynomial $y=a(x-b)^N$ used to fit the three points by least square method. In the first instance a uniform value of N=1 is used throughout the iteration which means the polynomial is a straight line which is a least square line fitted along the three equi-spaced points x, x+$\delta$ and x-$\delta$. In fact the convergence of the method using N=1 is very similar to the Newton method as will be shown in the examples provided.

In the second instance of the application of the three-point least square method, the value of N is allowed to dynamically vary with each step of the iteration. This procedure provides for additional flexibility since a better curve than straight line can be used as defined by the three points. Allowing N to vary with each step of the iteration is helpful in the initial steps of the iteration particularly for functions with higher gradients (derivatives). Towards the end of the iteration the value of N converges to N=1 in all cases.

In order to avoid off-shooting away from a possible nearer root by the use of too high value of N at any step of the iteration, it is possible to limit the variation of N to within the range:

-3 ≤ N ≤ 3

The stopping criterion used for the iteration process is given by:

$$|x_k - x_{k-1}| + |y_k| < 10^{-15} \qquad (61)$$

The rate of convergence towards the root x = r for each step of the iteration is evaluated using the formula:

$$C_k = \frac{Log\,(|E_{k+1}|)}{Log\,(|E_k|)} = \frac{Log\,(|X_{k+1} - r|)}{Log\,(|X_k - r|)} \qquad (62)$$

A quadratic convergence proved for this method is mostly evident with a $C_k$ value being close to 2 during the iteration.

The results of the iteration towards the root for seven equations shown in Table 1 are summarized along with the results of the use of Newton and secant methods for the purpose of comparison. Figure 3 shows an graphical display of the number of iterations required for the different equations tested.

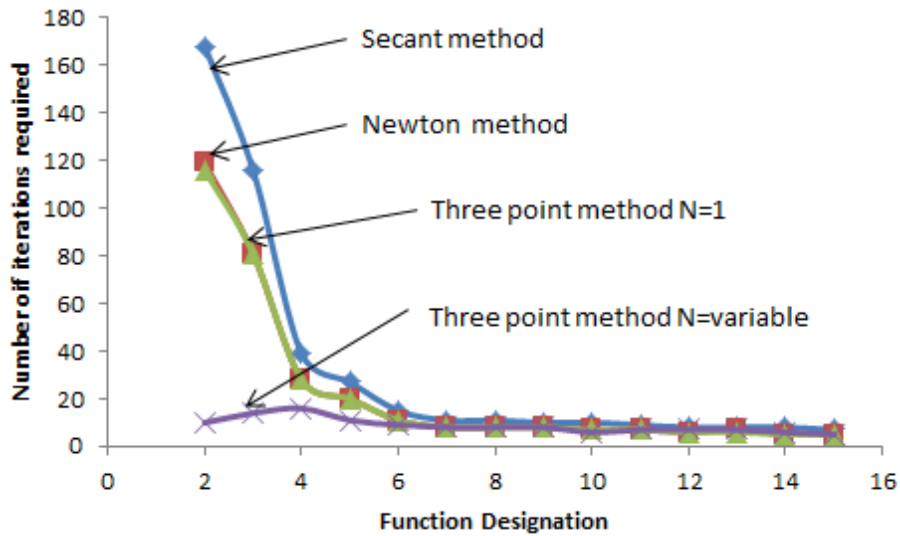

Figure 3  Comparison of the number of iterations required for the different equations tested.

Referring to table 1, the number of iterations required for the proposed method is equal to or less than that of Newton's method. For N=1 , the number of iterations required are more or less the same as that of Newton's method in almost all equations tested.  For the variable N case, better advantaged is provided for functions with higher gradients such as $y = (x − 2)(x + 2)^4$ and $y = e^{x^2+7x-30} − 1$ as the number of iterations required is significantly reduced. The secant method, having a less than quadratic rate of convergence, required in most cases the greatest number of iterations.

Table 1. Comparison of result of iteration of the three point method with Newton and secant methods.

| Function | Root | Starting point | Comparison of number of iterations required | | | |
|---|---|---|---|---|---|---|
| | | | secant Method | Newton Method | Least Square 3 –point Method | Least square 3 –point Method |
| | | | | | N= 1 | N = Variable |
| $y = x^3 + 4x^2 − 10$ | 1.365230013414100 | 0.5 | 10 | 8 | 8 | 8 |
| | | 1 | 8 | 6 | 6 | 7 |
| $y = [\sin(x)]^2 − x^2 + 1$ | -1.404491648215340 | -1 | 9 | 7 | 7 | 7 |
| | | -3 | 10 | 7 | 7 | 6 |
| $y = (x − 2)(x + 2)^4$ | -2.0000000000000 | -3 | 168 | 119 | 116 | 10 |
| | | 1.4 | 116 | 81 | 81 | 14 |
| | | 1.5 | 252 | 16 | 15 | 10 |
| $y = (x − 1)^6 − 1$ | 2.00000000000000 | 2.5 | 11 | 8 | 8 | 8 |
| | | 3.5 | 15 | 11 | 11 | 9 |
| $y = \sin(x).e^x + \ln(x^2 + 1)$ | -0.603231971557215 | -0.8 | 8 | 7 | 6 | 7 |
| | | -0.65 | 8 | 5 | 5 | 6 |
| $y = e^{x^2+7x-30} − 1$ | 3.000000000000000 | 4 | 27 | 20 | 20 | 11 |
| | | 4.5 | 39 | 28 | 28 | 16 |
| $y = x − 3 \ln(x)$ | 1.857183860207840 | 2 | 7 | 5 | 5 | 5 |
| | | 0.5 | 11 | 8 | 8 | 8 |

**Examples where the proposed method works while Newton method fails**

The advantage of the use of variable N is best illustrated by the application of the method where the Newton method and in several cases also the secant method fail to converge to the root.  The failure could be due to oscillation, divergence or off shooting to an invalid domain.  Table 2 below shows the results of the iterative process for the given equations where the newton and secant methods fail to converge with the starting points also indicated in the table.

As shown in Table 2, the proposed three point method with variable N does not result in failure to converge in all cases whereas the same method with N=1 shows failure in most of the cases where the Newton method fails also. This illustrates the advantage of using variable N rather than using N=1 for such non-convergent cases. This result also illustrates how the proposed method with fixed N (N=1) is closely similar to Newton's method in terms of both failure as well as rate of convergence.

Table 2. Results of application of the method for cases Newton or secant method fail to converge to the root.

| Function | Root | Starting point | secant Method | Newton Method | Least Square 3 –point Method N= 1 | Least square 3 – point Method N = Variable |
|---|---|---|---|---|---|---|
| $y = 2x^5 - 3x^4 + 4x^3 - x^2 + 10x - 13$ | 1.053392031515730 | 3.0 | 13 | Oscillates | 10 | 7 |
|  |  | -2.5 | 14 | Oscillates | 11 | 8 |
| $y = \log(x)$ | 1.000000000000000 | 3.0 | Fails | Fails | Fails | 7 |
| $y = Arctan(x)$ | 0.0000000000000 | 3.0 | Diverges | Diverges | Diverges | 7 |
|  |  | -3.0 | Diverges | Diverges | Diverges | 7 |
| $y = x^5 - x + 1$ | -1.167303978261420 | 2.0 | 48 | Oscillates | Oscillates | 10 |
|  |  | -3.0 | 14 | Oscillates | 11 | 7 |
| $y = 0.5x^3 - 6x^2 + 21.5x - 22$ | 4.00000000000000 | 3.0 | 7 | Oscillates | Oscillates | 7 |
| $y = x^{1/3}$ | 0.00000000000000 | 1.0 | Oscillates | Diverges | Diverges | 14 |
|  |  | -1.0 | Oscillates | Diverges | Diverges | 14 |
| $y = 10xe^{-x^2} - 1$ | 1.679630610428450 | 3.0 | Diverges | Diverges | Diverges | 11 |
|  | 0.101025848315685 | -1.0 | Diverges | Diverges | Diverges | 13 |

**Limitation of the proposed method**

In some cases during the iteration it might appear that y(x+δ) = y(x-δ). In this case because y(x+δ) - y(x-δ)=0, this results in division by zero and the value of δ should be readjusted to avoid such cases. However, this will not halt the iteration but calls for readjusting the value of δ such that y(x+δ) ≠y(x-δ).

The method requires evaluation of function for three points in each step of the iteration. In this regard the number of function evaluations required per each step of the iteration is higher than Newton and secant methods.

## 4. Conclusion

A method of root finding has been presented using a numerical iterative process involving three points together with a discussion of the derivation and proof of quadratic convergence. A given function $Y = f(x)$ whose root(s) are desired is fitted and approximated by a polynomial function curve of the form $y = a(x-b)^N$ and passing through three equi-spaced points using the principle of least squares. The method does not require evaluation of derivatives and requires only functional evaluations. The method has a quadratic convergence. The power of the polynomial curve used to fit the three equi-spaced points by least square method can be dynamically varied at each step of the iteration in order to provide better convergence characteristics or avoid oscillation, divergence and off shooting out of the valid domain for functional evaluation. From functional evaluation of the three equi-spaced points it is possible to make an estimate of the power N beforehand to be used in the next step of the iteration. An alternative application of the method using a uniform power of N=1 also gives a satisfactory result in many cases.

The limitation of the method is the necessity to evaluate the function at three points within each step of the iteration and the need to guard and alter the value of interval $\delta$ such that division by zero is avoided in the event $y(x+\delta) = y(x-\delta)$. However, this will not halt the iterative process only requiring adjusting the $\delta$ value.